\documentclass[12pt,intlimits,reqno]{article}

\usepackage{amssymb}
\usepackage{epsfig}
\usepackage{graphicx}
\usepackage{lscape}
\usepackage{amsmath}
\usepackage{amsfonts}
\usepackage{amscd}
\usepackage{cite}
\usepackage{amsthm}

\usepackage[all,cmtip]{xy}

\def\figurename{Figure}
\makeatletter
\renewcommand{\fnum@figure}[1]{\figurename~\thefigure.}

\newcommand{\RRe}{\operatorname{Re}}

\newtheorem{lemma}{\fontsize{14}{16pt}\selectfont Lemma}[section]
\newtheorem{theorem}{\fontsize{14}{16pt}\selectfont Theorem}[section]

\newtheorem*{note}{\fontsize{14}{16pt}\selectfont Note}
\newtheorem*{hypothesis}{\fontsize{14}{16pt}\selectfont Hypothesis}

\begin{document}

\thispagestyle{empty}

\begin{center}
    {\huge \textbf{New Sense of a Circle}} \\[1.5cm]
    {\Large \textbf{Mamuka Meskhishvili}}
\end{center}

\vskip+1.5cm

\begin{center}
    {\fontsize{14}{16pt}\selectfont \textbf{Abstract} }
\end{center}
\vskip+0.5cm

New condition is found for the set of points in the plane, for which the locus is a circle.

It is proved: the locus of points, such that the sum of the {\fontsize{14}{16pt}\selectfont $(2m)$}-th powers {\fontsize{14}{16pt}\selectfont $S_n^{(2m)}$} of the distances to the vertexes of fixed regular {\fontsize{14}{16pt}\selectfont $n$}-sided polygon is constant, is a circle if
\vskip+0.04cm
{\fontsize{14}{16pt}\selectfont
$$  S_n^{(2m)}>nr^{2m}, \;\;\text{\fontsize{12}{14pt}\selectfont where}\;\; m=1,2,\dots,n-1           $$
 }
\vskip+0.02cm
\noindent
and {\fontsize{14}{16pt}\selectfont $r$} is the distance from the center of the regular polygon to the vertex.

The radius {\fontsize{14}{16pt}\selectfont $\ell$} satisfies:
\vskip+0.02cm
{\fontsize{14}{16pt}\selectfont
$$  S_n^{(2m)}=n\Bigg[(r^2+\ell^2)^m+\sum_{k=1}^{[\frac{m}{2}]}
            \left(\!\!\!\begin{array}{c} \text{\fontsize{13}{14pt}\selectfont $m$} \\ \text{\fontsize{13}{14pt}\selectfont $2k$} \end{array}\!\!\!\right)
                         (r^2+\ell^2)^{m-2k}(r\ell)^{2k}
            \left(\!\!\!\begin{array}{c} \text{\fontsize{13}{14pt}\selectfont $2k$} \\ \text{\fontsize{13}{14pt}\selectfont $k$} \end{array}\!\!\!\right)\Bigg].       $$
 }

\vskip+1.5cm


\vskip+0.7cm

\textbf{2010 AMS Classification.}
    51M04, 51N20, 51N35.

\newpage

\ 

\bigskip
\bigskip

\noindent {\fontsize{18}{20pt}\selectfont \textbf{Contents}}

\vskip+0.7cm

\begin{description}

\item[1] Introduction \dotfill{3}

\item[2] Regular Polygon. Constant Sum of $2$nd Powers \dotfill{4}

\item[3] Regular Triangle. Constant Sum of $4$th Powers \dotfill{13}

\item[4] Square. Constant Sum of $4$th and $6$th Powers \dotfill{15}

\item[5] Regular Pentagon. Constant Sum of $4$th, $6$th and $8$th Powers \dotfill{18}

\item[6] Circle as Locus of Constant Sum of Powers \dotfill{26}

\item[7] Cosine Multiple Arguments and Cosine Powers Sums \dotfill{30}

\item[8] General Theorem \dotfill{34}
\end{description}

\newpage

\section{\fontsize{18}{20pt}\selectfont \textbf{Introduction}}
\vskip+1cm

A locus of points is the set of points, and only those points, that satisfies given conditions. A circle is usually defined as the locus of points (in the plane) at a given distance from a given point. This is well-known definition of a circle, but not only. It is known definition of a circle by using two fixed points -- the locus of points, such that the sum of the squares of the distances to two distinct fixed points is constant, is a circle.

More generally, for any collection of {\fontsize{14}{16pt}\selectfont $n$} number of points {\fontsize{14}{16pt}\selectfont $P_i$}, constants {\fontsize{14}{16pt}\selectfont $\lambda_i$} and {\fontsize{14}{16pt}\selectfont $S_n^{(2)}$}, distances {\fontsize{14}{16pt}\selectfont $d(X,P_i)$} -- the locus of points {\fontsize{14}{16pt}\selectfont $X$} such that
\vskip+0.02cm
{\fontsize{14}{16pt}\selectfont
$$  \sum_{i=1}^n \lambda_i d^2(X,P_i)=S_n^{(2)},    $$
 }
\vskip+0.2cm
\noindent is either
\begin{itemize}
\item[\textsf{(a)}] a circle, a point or empty set if {\fontsize{14}{16pt}\selectfont $\displaystyle \sum\limits_{i=1}^n \lambda_i\neq 0$},

\item[\textsf{(b)}] a line, a plane or the empty set if {\fontsize{14}{16pt}\selectfont $\displaystyle \sum\limits_{i=1}^n \lambda_i=0$}.
\end{itemize}

Indeed, let {\fontsize{14}{16pt}\selectfont $(a_i,b_i)$} be the coordinates of point {\fontsize{14}{16pt}\selectfont $P_i$} and {\fontsize{14}{16pt}\selectfont $(x,y)$} the coordinates of point {\fontsize{14}{16pt}\selectfont $X$}. Then the equation satisfied by {\fontsize{14}{16pt}\selectfont $X$} takes the form
\vskip+0.02cm
{\fontsize{14}{16pt}\selectfont
\begin{align*}
    S_n^{(2)} & =\sum_{i=1}^n \lambda_i\big((x-a_i)^2+(y-b_i)^2\big)= \\[0.2cm]
    & =(x^2+y^2)\sum_{i=1}^n \lambda_i -x\cdot 2\sum_{i=1}^n \lambda_ia_i- \\
    &\qquad\qquad -y\cdot 2\sum_{i=1}^n \lambda_ib_i+\sum_{i=1}^n\lambda_i (a_i^2+b_i^2).
\end{align*}
 }
\vskip+0.2cm
\noindent
If the coefficient of {\fontsize{14}{16pt}\selectfont $x^2+y^2$} is nonzero, then this equation determines either a circle, a point or the empty set and if it is zero, then the equation determines either a line, a plane, or the empty set {\fontsize{14}{16pt}\selectfont [1]}.

In this article we consider under points {\fontsize{14}{16pt}\selectfont $P_i$} -- vertexes of a regular {\fontsize{14}{16pt}\selectfont $n$}-sided polygon {\fontsize{14}{16pt}\selectfont $P_n$} and for coefficients
\vskip+0.02cm
{\fontsize{14}{16pt}\selectfont
$$  \lambda_i=1, \;\; i=1,2,\dots,n.      $$
 }

We start {\fontsize{14}{16pt}\selectfont $m=1$}, i.e. we have the constant sum of the second powers of distances from point {\fontsize{14}{16pt}\selectfont $X$} to the vertexes of regular {\fontsize{14}{16pt}\selectfont $n$}-sided polygon. In our notations -- case {\fontsize{14}{16pt}\selectfont $S_n^{(2)}$}. The case of two points {\fontsize{14}{16pt}\selectfont ($n=2$)} is well-known {\fontsize{14}{16pt}\selectfont [2]}, thus we discuss {\fontsize{14}{16pt}\selectfont $n>2$}.

\vskip+1.5cm
\section{\fontsize{18}{20pt}\selectfont \textbf{Regular Polygon. \\ Constant Sum of 2nd Powers}}
\vskip+1cm


\begin{theorem}\label{th:1}
{\fontsize{14}{16pt}\selectfont The locus of points, such that the sum of the  \linebreak   se\-cond powers {\fontsize{14}{16pt}\selectfont $S_n^{(2)}$} of the distances to the vertexes of fixed regular  \linebreak   {\fontsize{14}{16pt}\selectfont $n$}-sided polygon is constant, is a circle if
\vskip+0.02cm
{\fontsize{14}{16pt}\selectfont
$$  S_n^{(2)}>nr^2,     $$
 }
\vskip+0.2cm
\noindent
where {\fontsize{14}{16pt}\selectfont $r$} is the distance from the center of the regular polygon to the vertex. The center of the circle is at the center of the regular polygon and the radius is
\vskip+0.02cm
{\fontsize{14}{16pt}\selectfont
$$  \sqrt{\frac{S_n^{(2)}}{n}-r^2}\,.           $$
 }
\vskip+0.2cm

If {\fontsize{14}{16pt}\selectfont $S_n^{(2)}=nr^2$} the locus is a point -- the center of polygon and if {\fontsize{14}{16pt}\selectfont $S_n^{(2)}<nr^2$} it is the empty set. }
\end{theorem}
\bigskip

\begin{figure}[h]
\centerline{\includegraphics[width=7.5cm]
    {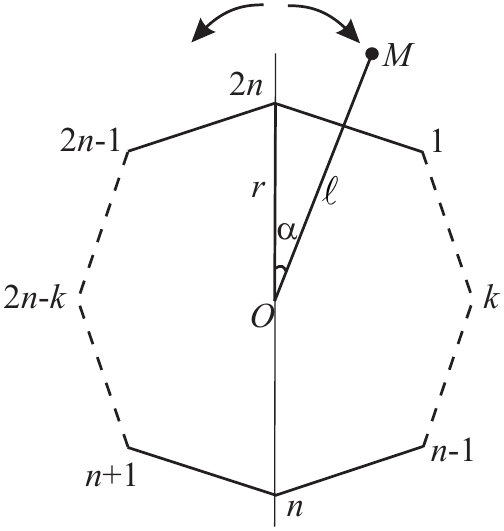}}
\caption{}
\label{fig:1}
\end{figure}

Firstly discuss the regular polygon with even number of vertexes {\fontsize{14}{16pt}\selectfont $2n$}. Numerate vertexes as shown in Figure \ref{fig:1}. Denote the distance from an arbitrary point {\fontsize{14}{16pt}\selectfont $M$} to the vertex {\fontsize{14}{16pt}\selectfont $i$} by {\fontsize{14}{16pt}\selectfont $x_i$}, respectively, {\fontsize{14}{16pt}\selectfont $i=1,\dots,2n$} and {\fontsize{14}{16pt}\selectfont $OM=\ell$}. In even case vertexes
\vskip+0.2cm
{\fontsize{14}{16pt}\selectfont
$$  1 \;\text{\fontsize{12}{14pt}\selectfont and}\; 2n-1\;\;,\dots,\;\;k \;\text{\fontsize{12}{14pt}\selectfont and}\; 2n-k\;\;,\dots,\;\;n-1 \; \text{\fontsize{12}{14pt}\selectfont and}\;
         n+1        $$
 }
  \vskip+0.2cm
 \noindent are symmetrical with respect to the line which passes through center {\fontsize{14}{16pt}\selectfont $0$} and vertex {\fontsize{14}{16pt}\selectfont $2n$}, as shown in Figure \ref{fig:1}.

\noindent

\pagebreak
Calculate the square of distances in pairs
\vskip+0.02cm
{\fontsize{14}{16pt}\selectfont
\begin{gather*}
    x_1^2 \;\text{\fontsize{12}{14pt}\selectfont and}\; x_{2n-1}^2 \;\;,\dots, \;\; x_k^2 \;\text{\fontsize{12}{14pt}\selectfont and}\; x_{2n-k}^2 \;\;,\dots, \;\; x_{n-1}^2 \;\text{\fontsize{12}{14pt}\selectfont and}\; x_{n+1}^2, \\[0.2cm]
    x_n^2 \;\text{\fontsize{12}{14pt}\selectfont and}\; x_{2n}^2.
\end{gather*}
 }
\vskip+0.2cm
\noindent
To calculate angles, we use arrows directions (Figure \ref{fig:1}) clockwise for    \linebreak        {\fontsize{14}{16pt}\selectfont $1,\dots,n-1$} vertexes and counterclockwise for {\fontsize{14}{16pt}\selectfont $2n-1,\dots,n+1$} vertexes;
\vskip+0.02cm
{\fontsize{14}{16pt}\selectfont
\allowdisplaybreaks
\begin{align*}
    x_1^2& =r^2+\ell^2-2r\ell\cos\Big(\frac{360^\circ}{2n}-\alpha\Big) \;\;\text{\fontsize{12}{14pt}\selectfont and}\;\; \\[0.2cm]
    &\qquad\qquad\qquad x_{2n-1}^2=r^2+\ell^2-2r\ell\cos\Big(\frac{360^\circ}{2n}+\alpha\Big), \\[0.2cm]
    &\qquad\qquad \vdots \\[0.2cm]
    x_k^2 & =r^2+\ell^2-2r\ell\cos\Big(k\,\frac{360^\circ}{2n}-\alpha\Big) \;\;\text{\fontsize{12}{14pt}\selectfont and}\;\; \\[0.2cm]
    &\qquad\qquad\qquad x_{2n-k}^2=r^2+\ell^2-2r\ell\cos\Big(k\,\frac{360^\circ}{2n}+\alpha\Big), \\[0.2cm]
    &\qquad\qquad \vdots \\[0.2cm]
    x_{n-1}^2 & =r^2+\ell^2-2r\ell\cos\Big((n-1)\,\frac{360^\circ}{2n}-\alpha\Big) \;\;\text{\fontsize{12}{14pt}\selectfont and}\;\; \\[0.2cm]
    &\qquad\qquad\qquad x_{n+1}^2=r^2+\ell^2-2r\ell\cos\Big((n-1)\,\frac{360^\circ}{2n}+\alpha\Big), \\[0.2cm]
    x_n^2 & =r^2+\ell^2-2r\ell\cos(180^\circ-\alpha) \;\;\text{\fontsize{12}{14pt}\selectfont and}\;\; \\[0.2cm]
    &\qquad\qquad\qquad x_{2n}^2=r^2+\ell^2-2r\ell\cos\alpha.
\end{align*}
 }
\vskip+0.2cm
\noindent
Summarize each pair:
\vskip+0.02cm
{\fontsize{14}{16pt}\selectfont
\allowdisplaybreaks
\begin{align*}
    x_1^2+x_{2n-1}^2 & =2(r^2+\ell^2)-2r\ell\cdot 2\cos\alpha\cdot\cos\frac{360^\circ}{2n}\,, \\[0.2cm]
    &\qquad \vdots \\[0.2cm]
    x_k^2+x_{2n-k}^2 & =2(r^2+\ell^2)-2r\ell\cdot 2\cos\alpha\cdot\cos k\,\frac{360^\circ}{2n}\,, \\[0.2cm]
    &\qquad \vdots \\[0.2cm]
    x_{n-1}^2+x_{n+1}^2 & =2(r^2+\ell^2)-2r\ell\cdot 2\cos\alpha\cdot\cos\Big((n-1)\,\frac{360^\circ}{2n}\Big), \\[0.2cm]
    x_n^2+x_{2n}^2 & =2(r^2+\ell^2)-2r\ell\cdot 2\cos 90^\circ\cdot\cos(90^\circ-\alpha)= \\[0.2cm]
    & =2(r^2+\ell^2).
\end{align*}
 }
\vskip+0.2cm
\noindent
There are {\fontsize{14}{16pt}\selectfont $n$} pairs, so total sum is
\vskip+0.02cm
{\fontsize{14}{16pt}\selectfont
\allowdisplaybreaks
\begin{align*}
    S_{2n}^{(2)} & =\sum_{i=1}^{2n} x_i^2= \\[0.2cm]
    & =2n(r^2+\ell^2)-4r\ell\cos\alpha\sum_{k=1}^{n-1} \cos\Big(k\,\frac{360^\circ}{2n}\Big).
\end{align*}
\begin{multline*}
    \sum_{k=1}^{n-1} \cos\Big(k\,\frac{360^\circ}{2n}\Big)=\sum_{k=1}^{n-1} \cos\Big(k\,\frac{180^\circ}{n}\Big)= \\[0.2cm]
    =\cos\frac{180^\circ}{n}+\cos 2\,\frac{180^\circ}{n}+\cdots+\cos (n-2)\,\frac{180^\circ}{n}+\cos (n-1)\,\frac{180^\circ}{n}\,.
\end{multline*}
 }
\vskip+0.2cm
\noindent
After grouping
\vskip+0.02cm
{\fontsize{14}{16pt}\selectfont
\allowdisplaybreaks
\begin{multline*}
    \sum_{k=1}^{n-1} \cos\Big(k\,\frac{180^\circ}{n}\Big)= \\[0.2cm]
    =\cos\frac{180^\circ}{n}+\cos (n-1)\,\frac{180^\circ}{n}+\cos 2\,\frac{180^\circ}{n}+
                \cos (n-2)\,\frac{180^\circ}{n}+\cdots+ \\[0.2cm]
    =\cos\frac{180^\circ}{n}+\cos\Big(180^\circ-\frac{180^\circ}{n}\Big)+\cos 2\,\frac{180^\circ}{n}+ \\[0.2cm]
    +\cos\Big(180^\circ-2\,\frac{180^\circ}{n}\Big)+\cdots=0,
\end{multline*}
 }
\vskip+0.2cm
\noindent
since each pair is {\fontsize{14}{16pt}\selectfont $0$}.

\bigskip
\noindent \textbf{Note.} If {\fontsize{14}{16pt}\selectfont $n$} is odd we have {\fontsize{14}{16pt}\selectfont $\dfrac{n-1}{2}$} pairs, but if {\fontsize{14}{16pt}\selectfont $n$} is even {\fontsize{14}{16pt}\selectfont $\dfrac{n}{2}-1$} pairs (whose sums are {\fontsize{14}{16pt}\selectfont $0$}) and one more alone term {\fontsize{14}{16pt}\selectfont $\displaystyle \cos\Big(\frac{n}{2}\cdot\frac{180^\circ}{n}\Big)$}\,, which itself is {\fontsize{14}{16pt}\selectfont $0$}. So we prove
\vskip+0.02cm
{\fontsize{14}{16pt}\selectfont
$$  \sum_{k=1}^{n-1} \cos k\,\frac{180^\circ}{n}=0,     $$
 }
\vskip+0.2cm
\noindent
and from that it follows
\vskip+0.02cm
{\fontsize{14}{16pt}\selectfont
$$  S_n^{(2)}=\sum_{i=1}^{2n} x_i^2=2n(r^2+\ell^2).     $$
 }
\vskip+0.2cm

Now discuss the case, when the number of vertexes of a regular polygon is odd -- {\fontsize{14}{16pt}\selectfont $2n+1$}.

\begin{figure}[h]
\centerline{\includegraphics[width=7.5cm]{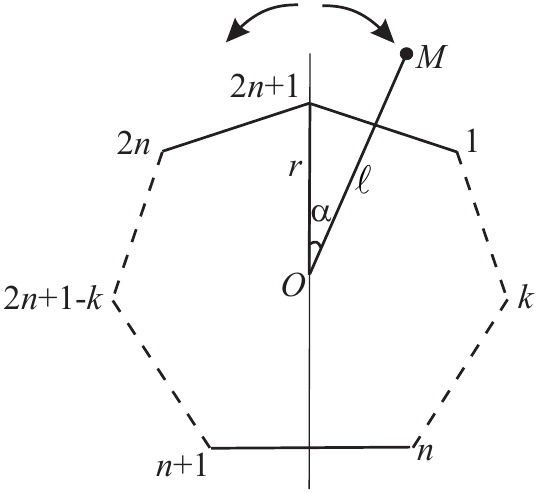}}
\caption{}
\label{fig:2}
\end{figure}

Numerate vertexes as shown in Figure \ref{fig:2}. As in even case vertexes
\vskip+0.02cm
{\fontsize{14}{16pt}\selectfont
$$  1 \;\text{\fontsize{12}{14pt}\selectfont and}\; 2n \;\;,\dots,\;\; k \;\text{\fontsize{12}{14pt}\selectfont and}\; 2n+1-k \;\;,\dots, \;\; n
 \;\text{\fontsize{12}{14pt}\selectfont and}\; n+1      $$
 }
 \vskip+0.2cm
\noindent
are symmetrical with respect to the line which passes through center {\fontsize{14}{16pt}\selectfont $0$} and vertex {\fontsize{14}{16pt}\selectfont $2n+1$}, as shown in Figure \ref{fig:2}.

\noindent Calculate the square of distances in pairs --
\vskip+0.02cm
{\fontsize{14}{16pt}\selectfont
$$  x_1^2 \;\text{\fontsize{12}{14pt}\selectfont and}\; x_{2n}^2 \;\;,\dots,\;\; x_k^2 \;\text{\fontsize{12}{14pt}\selectfont and}\; x_{2n+1-k}^2 \;\;,\dots,\;\;
            x_n^2 \;\text{\fontsize{12}{14pt}\selectfont and}\; x_{n+1}^2.   $$
 }
\vskip+0.2cm
\noindent
Only {\fontsize{14}{16pt}\selectfont $x_{2n+1}^2$} -- square of distance remains single. To calculate the angles, we use arrows directions (Figure \ref{fig:2}) clockwise for vertexes -- {\fontsize{14}{16pt}\selectfont $1,\dots,n$} and counterclockwise for vertexes {\fontsize{14}{16pt}\selectfont $2n,\dots,n+1$};
\vskip+0.02cm
{\fontsize{14}{16pt}\selectfont
\allowdisplaybreaks
\begin{align*}
    x_1^2 & =r^2+\ell^2-2r\ell\cos\Big(\frac{360^\circ}{2n+1}-\alpha\Big) \;\;\text{\fontsize{12}{14pt}\selectfont and}\;\; \\[0.2cm]
    &\qquad\qquad\qquad x_{2n}^2=r^2+\ell^2-2r\ell\cos\Big(\frac{360^\circ}{2n+1}+\alpha\Big), \\[0.2cm]
    &\qquad\qquad \vdots \\[0.2cm]
    x_k^2 & =r^2+\ell^2-2r\ell\cos\Big(k\,\frac{360^\circ}{2n+1}-\alpha\Big) \;\;\text{\fontsize{12}{14pt}\selectfont and}\;\; \\[0.2cm]
    &\qquad\qquad\qquad x_{2n+1-k}^2=r^2+\ell^2-2r\ell\cos\Big(k\,\frac{360^\circ}{2n+1}+\alpha\Big), \\[0.2cm]
    &\qquad\qquad \vdots \\[0.2cm]
    x_n^2 & =r^2+\ell^2-2r\ell\cos\Big(n\,\frac{360^\circ}{2n+1}-\alpha\Big) \;\;\text{\fontsize{12}{14pt}\selectfont and}\;\; \\[0.2cm]
    &\qquad\qquad\qquad x_{n+1}^2=r^2+\ell^2-2r\ell\cos\Big(n\,\frac{360^\circ}{2n+1}+\alpha\Big), \\[0.2cm]
    x_{2n+1}^2 & =r^2+\ell^2-2r\ell\cos\alpha.
\end{align*}
 }
\vskip+0.2cm
\noindent
There are {\fontsize{14}{16pt}\selectfont $n$} pairs, so total sum is
\vskip+0.02cm
{\fontsize{14}{16pt}\selectfont
\allowdisplaybreaks
\begin{align*}
    S_{2n+1}^{(2)} & =\sum_{i=1}^{2n+1} x_i^2= \\[0.2cm]
    & =(2n+1)(r^2+\ell^2)-2r\ell\cos\alpha- \\[0.2cm]
    &\quad -2r\ell\sum_{k=1}^n \bigg(\cos\Big(k\,\frac{360^\circ}{2n+1}-\alpha\Big)+\cos\Big(k\,\frac{360^\circ}{2n+1}+\alpha\Big)\bigg)= \\[0.2cm]
    & =(2n+1)(r^2+\ell^2)-2r\ell\cos\alpha- \\[0.2cm]
    &\qquad -2r\ell\sum_{k=1}^n 2\cos\alpha\cdot\cos\Big(k\,\frac{360^\circ}{2n+1}\Big)= \\[0.2cm]
    & =(2n+1)(r^2+\ell^2)- \\[0.2cm]
    &\qquad -2r\ell\cos\alpha\bigg(1+2\sum_{k=1}^n \Big(\cos k\,\frac{360^\circ}{2n+1}\Big)\bigg).
\end{align*}
 }
\vskip+0.2cm
\noindent
Each addend of the sum
\vskip+0.02cm
{\fontsize{14}{16pt}\selectfont
$$  \sum_{k=1}^n \cos k\,\frac{360^\circ}{2n+1}\,,      $$
 }
\vskip+0.2cm
\noindent
multiply by
\vskip+0.02cm
{\fontsize{14}{16pt}\selectfont
$$   2\sin\frac{180^\circ}{2n+1}  \,,      $$
}
\vskip+0.2cm
\noindent
and then use trigonometric identity
\vskip+0.02cm
{\fontsize{14}{16pt}\selectfont
$$  2\sin C\cdot \cos D=\sin(C+D)-\sin(D-C)     $$
 }
\vskip+0.2cm
\noindent
we get
\vskip+0.02cm
{\fontsize{14}{16pt}\selectfont
\allowdisplaybreaks
\begin{align*}
    2\sin & \,\frac{180^\circ}{2n+1}\,\Big(\cos\frac{2\cdot 180^\circ}{2n+1}+\cos 2\,\frac{2\cdot 180^\circ}{2n+1}+  \\[0.5cm]
    &\qquad\qquad\qquad\qquad + \cos 3\,\frac{2\cdot 180^\circ}{2n+1}+\cdots+ \\[0.5cm]
    &\qquad\qquad +\cos(n-1)\,\frac{2\cdot 180^\circ}{2n+1}+\cos n\,\frac{2\cdot 180^\circ}{2n+1}\Big)= 
\end{align*}
\newpage
\begin{align*}
    & =\sin\frac{3\cdot 180^\circ}{2n+1}-\sin\frac{180^\circ}{2n+1}+ \\[0.2cm]
    &\qquad +\sin\frac{5\cdot 180^\circ}{2n+1}-\sin\frac{3\cdot 180^\circ}{2n+1}+ \\[0.2cm]
    &\qquad +\sin\frac{7\cdot 180^\circ}{2n+1}-\sin\frac{5\cdot 180^\circ}{2n+1}+ \\[0.2cm]
    &\qquad\qquad\qquad \vdots \\[0.2cm]
    &\qquad +\sin\frac{(2n-1)\cdot 180^\circ}{2n+1}-\sin\frac{(2n-3)\cdot 180^\circ}{2n+1}+ \\[0.2cm]
    &\qquad +\sin\frac{(2n+1)\cdot 180^\circ}{2n+1}-\sin\frac{(2n-1)\cdot 180^\circ}{2n+1}= \\[0.2cm]
    & =\sin\frac{(2n+1)\cdot 180^\circ}{2n+1}-\sin\frac{180^\circ}{2n+1}=-\sin\frac{180^\circ}{2n+1}\,.
\end{align*}
 }
\vskip+0.2cm
\noindent
So,
{\fontsize{14}{16pt}\selectfont
$$  \sum_{k=1}^n \cos k\,\frac{360^\circ}{2n+1}=-\frac{1}{2}        $$
 }
\vskip+0.2cm
\noindent
and in odd case we obtain
\vskip+0.02cm
{\fontsize{14}{16pt}\selectfont
$$  S_{2n+1}^{(2)}=\sum_{i=1}^{2n+1} x_i^2=(2n+1)(r^2+\ell^2).      $$
 }
\vskip+0.2cm

Summarize even and odd cases. For an arbitrary number {\fontsize{14}{16pt}\selectfont $n$}, we have:
\vskip+0.02cm
{\fontsize{14}{16pt}\selectfont
$$  S_n^{(2)}=n(r^2+\ell^2),       $$
 }
 \vskip+0.2cm
\noindent 
which finally proves the Theorem \ref{th:1}.

\pagebreak
\section{\fontsize{18}{20pt}\selectfont \textbf{Regular Triangle. \\ Constant Sum of 4th Powers}}
\vskip+1cm

The sum of the fourth powers of the distances to the vertexes is:
\vskip+0.02cm
{\fontsize{14}{16pt}\selectfont
\begin{align*}
    S_3^{(4)} & =\sum_{i=1}^3 x_i^4= \\[0.2cm]
    & =(r^2+\ell^2-2r\ell\cos\alpha)^2+\big(r^2+\ell^2-2r\ell\cos(120^\circ-\alpha)\big)^2+ \\[0.2cm]
    &\qquad\qquad +\big(r^2+\ell^2-2r\ell\cos(120^\circ+\alpha)\big)^2.
\end{align*}
 }
\vskip+0.2cm

Introduce new notations
\vskip+0.02cm
{\fontsize{14}{16pt}\selectfont
\allowdisplaybreaks
\begin{gather*}
    A=r^2+\ell^2 \;\;\text{\fontsize{12}{14pt}\selectfont and}\;\; B=2r\ell; \\[0.3cm]
\begin{aligned}
    S_3^{(4)} & =(A-B\cos\alpha)^2+\big(A-B\cos(120^\circ-\alpha)\big)^2+ \\[0.2cm]
    &\qquad\;\; +\big(A-B\cos(120^\circ+\alpha)\big)^2= \\[0.2cm]
    & =3A^2-2AB\big(\cos\alpha+\cos(120^\circ-\alpha)+\cos(120^\circ+\alpha)\big)+ \\[0.2cm]
    &\qquad\;\; +B^2\big(\cos^2\alpha+\cos^2(120^\circ-\alpha)+\cos^2(120^\circ+\alpha)\big);
\end{aligned} \\[0.5cm]
    \cos\alpha+\cos(120^\circ-\alpha)+\cos(120^\circ+\alpha)= \qquad\qquad\qquad\qquad \\[0.2cm]
    \qquad\qquad\qquad =\cos\alpha+2\cos 120^\circ \cos\alpha=0, \\[0.2cm]
    \cos^2\alpha+\cos^2(120^\circ-\alpha)+\cos^2(120^\circ+\alpha)= \qquad\qquad\qquad\qquad \\[0.2cm]
    \qquad\qquad =\frac{1}{2}\,\Big(3+\cos 2\alpha+\cos(240^\circ-2\alpha)+\cos (240^\circ+2\alpha)\Big)= \\[0.2cm]
    =\frac{3}{2}+\frac{1}{2}\,\Big(\cos 2\alpha+\cos(120^\circ+2\alpha)+\cos(120^\circ-2\alpha)\Big)=\frac{3}{2}\,.
\end{gather*}
 }
\vskip+0.2cm
\noindent
So,
\vskip+0.02cm
{\fontsize{14}{16pt}\selectfont
$$  S_3^{(4)}=3A^2+\frac{3}{2}\,B^2=3(r^4+\ell^4+4r^2\ell^2).       $$
 }
\vskip+0.2cm

The sum {\fontsize{14}{16pt}\selectfont $S_3^{(4)}$} does not contain {\fontsize{14}{16pt}\selectfont $\alpha$}, i.e. this expression does not depend on direction {\fontsize{14}{16pt}\selectfont $OM$} and only depends on length -- {\fontsize{14}{16pt}\selectfont $OM$}.

Since the triangle is fixed {\fontsize{14}{16pt}\selectfont $r$} is constant and this expression is increasing under {\fontsize{14}{16pt}\selectfont $\ell$} {\fontsize{14}{16pt}\selectfont $(\ell>0)$}, so obtained relation is one to one:
\vskip+0.02cm
{\fontsize{14}{16pt}\selectfont
$$  S_3^{(4)}=const \;\;\Longleftrightarrow\;\; \ell=const.        $$
 }
\vskip+0.2cm

We proved

\bigskip
\begin{theorem}\label{th:2}
{\fontsize{14}{16pt}\selectfont
The locus of points, such that the sum of the fourth powers {\fontsize{14}{16pt}\selectfont $S_3^{(4)}$} of the distances to the vertexes of fixed regular triangle is constant, is a circle if
\vskip+0.02cm
{\fontsize{14}{16pt}\selectfont
$$  S_3^{(4)}>3r^4,     $$
 }
\vskip+0.2cm
\noindent
where {\fontsize{14}{16pt}\selectfont $r$} is the distance from the center of the regular triangle to the vertex. The center of the circle is at the center of the regular triangle and the radius {\fontsize{14}{16pt}\selectfont $\ell$} satisfies the following condition:
\vskip+0.02cm
{\fontsize{14}{16pt}\selectfont
$$  S_3^{(4)}=3\Big[(r^2+\ell^2)^2+2(r\ell)^2\Big].     $$
 }
\vskip+0.2cm

If {\fontsize{14}{16pt}\selectfont $S_3^{(4)}=3r^4$} the locus is a point -- the center of the regular triangle and if {\fontsize{14}{16pt}\selectfont $S_3^{(4)}<3r^4$} it is the empty set. }
\end{theorem}
\bigskip

\noindent \textbf{Note.} The sum of the sixth powers -- {\fontsize{14}{16pt}\selectfont $S_3^{(6)}$} contains {\fontsize{14}{16pt}\selectfont $\alpha$}, so in this case the locus is not a circle.

\vskip+1.5cm
\section{\fontsize{18}{20pt}\selectfont \textbf{Square. \\ Constant Sum of 4th and 6th Powers}}
\vskip+0.3cm
{\fontsize{14}{16pt}\selectfont
\allowdisplaybreaks
\begin{align*}
    S_4^{(4)} & =(A-B\cos\alpha)^2+\big(A-B\cos(90^\circ-\alpha)\big)^2+ \\[0.2cm]
    &\qquad +\big(A-B\cos(180^\circ-\alpha)\big)^2+\big(A-B\cos(90^\circ+\alpha)\big)^2= \\[0.2cm]
    & =4A^2-2AB\big(\cos\alpha+\cos(90^\circ-\alpha)+ \\[0.2cm]
    &\qquad\qquad\qquad +\cos(180^\circ-\alpha)+\cos(90^\circ+\alpha)\big)+ \\[0.2cm]
    &\qquad +B^2\big(\cos^2\alpha+\cos^2(90^\circ-\alpha)+ \\[0.2cm]
    &\qquad\qquad\qquad +\cos^2(180^\circ-\alpha)+\cos^2(90^\circ+\alpha)\big);
\end{align*} \\[-0.7cm]
\begin{gather*}
    \cos\alpha+\cos(90^\circ-\alpha)+\cos(180^\circ-\alpha)+\cos(90^\circ+\alpha)= \qquad\qquad\qquad \\[0.2cm]
\qquad\qquad\qquad =2\cos 90^\circ\cdot \cos(90^\circ-\alpha)+2\cos 90^\circ\cdot\cos\alpha=0, \\[1cm]
    \cos^2\alpha+\cos^2(90^\circ-\alpha)+\cos^2(180^\circ-\alpha)+\cos^2(90^\circ+\alpha)= \qquad\qquad\qquad \\[0.2cm]
    =\frac{1}{2}\,\Big(4+\cos 2\alpha+\cos(180^\circ-2\alpha)+ \qquad\qquad\qquad\qquad\qquad \\[0.2cm]
    \qquad\qquad\qquad +\cos(360^\circ-2\alpha)+\cos(180^\circ+2\alpha)\Big)=2, \\[0.5cm]
    S_4^{(4)}=4A^2+2B^2=4(r^4+\ell^4+4r^2\ell^2).
\end{gather*}
 }

\begin{theorem}\label{th:3}
{\fontsize{14}{16pt}\selectfont
The locus of points, such that the sum of the fourth powers {\fontsize{14}{16pt}\selectfont $S_4^{(4)}$} of the distances to the vertexes of fixed square is constant, is a circle if
\vskip+0.02cm
{\fontsize{14}{16pt}\selectfont
$$  S_4^{(4)}>4r^4,         $$
 }
\vskip+0.1cm
\noindent where {\fontsize{14}{16pt}\selectfont $r$} is the distance from the center of the square to the vertex. The center of the circle is at the center of the square and the radius {\fontsize{14}{16pt}\selectfont $\ell$} satisfying the following condition:
\vskip+0.02cm
{\fontsize{14}{16pt}\selectfont
$$  S_4^{(4)}=4\Big[(r^2+\ell)^2+2(r\ell)^2\Big].     $$
 }
\vskip+0.05cm

If {\fontsize{14}{16pt}\selectfont $S_4^{(4)}=4r^4$} the locus is a point -- the center of the square and if {\fontsize{14}{16pt}\selectfont $S_4^{(4)}<4r^4$} it is the empty set. }
\end{theorem}

\medskip
The sum of sixth powers:
\vskip+0.02cm
{\fontsize{14}{16pt}\selectfont
\allowdisplaybreaks
\begin{align*}
    S_4^{(6)} & =(A-B\cos\alpha)^3+\big(A-B\cos(90^\circ-\alpha)\big)^3+ \\[0.2cm]
    &\qquad +\big(A-B\cos(180^\circ-\alpha)\big)^3+\big(A-B\cos(90^\circ+\alpha)\big)^3= \\[0.2cm]
    & =4A^3-3A^2B\Big(\cos\alpha+\cos(90^\circ-\alpha)+ \\[0.2cm]
    &\qquad\qquad\qquad\qquad +\cos(180^\circ-\alpha)+\cos(90^\circ+\alpha)\Big)+ \\[0.2cm]
    &\qquad +3AB^2\Big(\cos^2\alpha+\cos^2(90^\circ-\alpha)+ \\[0.2cm]
    &\qquad\qquad\qquad\qquad +\cos^2(180^\circ-\alpha)+\cos^2(90^\circ+\alpha)\Big)- \\[0.2cm]
    &\qquad -B^3\Big(\cos^3\alpha+\cos^3(90^\circ-\alpha)+ \\[0.2cm]
    &\qquad\qquad\qquad\qquad +\cos^3(180^\circ-\alpha)+\cos^3(90^\circ+\alpha)\Big),
\end{align*}
\begin{multline*}
    \cos^3\alpha+\cos^3(90^\circ-\alpha)+\cos^3(180^\circ-\alpha)+\cos^3(90^\circ+\alpha)= \\[0.2cm]
    =\cos^3\alpha+\cos^3(90^\circ-\alpha)-\cos^3\alpha-\cos^3(90^\circ-\alpha)=0.
\end{multline*}
 }
\vskip+0.2cm
\noindent
So,
\vskip+0.02cm
{\fontsize{14}{16pt}\selectfont
$$  S_4^{(6)}=4A^3+3AB^2\cdot 2=4(r^2+\ell^2)(r^4+\ell^4+8r^2\ell^2).       $$
 }
\vskip+0.2cm
\noindent
For fixed {\fontsize{14}{16pt}\selectfont $r$}, obtained expression is increasing under variable {\fontsize{14}{16pt}\selectfont $\ell$} {\fontsize{14}{16pt}\selectfont $(\ell>0)$}, so it is one to one relation between the sum of the sixth powers {\fontsize{14}{16pt}\selectfont $S_4^{(6)}$} and {\fontsize{14}{16pt}\selectfont $\ell$}:
\vskip+0.02cm
{\fontsize{14}{16pt}\selectfont
$$  S_4^{(6)}=const \;\;\Longleftrightarrow\;\; \ell=const.     $$
 }

\begin{theorem}\label{th:4}
{\fontsize{14}{16pt}\selectfont
The locus of points, such that the sum of the sixth powers {\fontsize{14}{16pt}\selectfont $S_4^{(6)}$} of the distances to the vertexes of fixed square is constant, is a circle if
\vskip+0.02cm
{\fontsize{14}{16pt}\selectfont
$$  S_4^{(6)}>4r^6,     $$
 }
\vskip+0.2cm
\noindent
where {\fontsize{14}{16pt}\selectfont $r$} is the distance from the center of the square to the vertex. The center of the circle is at the center of the square and the radius {\fontsize{14}{16pt}\selectfont $\ell$} satisfying the following condition:
\vskip+0.02cm
{\fontsize{14}{16pt}\selectfont
$$  S_4^{(6)}=4\Big[(r^2+\ell^2)^3+6(r^2+\ell^2)(r\ell)^2\Big].     $$
 }
\vskip+0.2cm

If {\fontsize{14}{16pt}\selectfont $S_4^{(6)}=4r^6$} the locus is a point -- the center of the square and if {\fontsize{14}{16pt}\selectfont $S_4^{(6)}<4r^6$} it is the empty set. }
\end{theorem}

\bigskip
\noindent \textbf{Note.} The sum of the eighth powers -- {\fontsize{14}{16pt}\selectfont $S_4^{(8)}$} contains {\fontsize{14}{16pt}\selectfont $\alpha$}, so in this case the locus is not a circle.

\vskip+1.5cm
\section{\fontsize{18}{20pt}\selectfont \textbf{Regular Pentagon. \\ Constant Sum of 4th, 6th and 8th Powers}}
\vskip+0.3cm
{\fontsize{14}{16pt}\selectfont
\allowdisplaybreaks
\begin{align*}
    S_5^{(4)} & =(r^2+\ell^2-2r\ell\cos\alpha)^2+\big(r^2+\ell^2-2r\ell\cos(72^\circ-\alpha)\big)^2+ \\[0.2cm]
    &\qquad\;\; +\big(r^2+\ell^2-2r\ell\cos(144^\circ-\alpha)\big)^2+ \\[0.2cm]
    &\qquad\;\; +\big(r^2+\ell^2-2r\ell\cos(144^\circ+\alpha)\big)^2+ \\[0.2cm]
    &\qquad\;\; +\big(r^2+\ell^2-2r\ell\cos(72^\circ+\alpha)\big)^2= \\[0.2cm]
    & =(A-B\cos\alpha)^2+\big(A-B\cos(72^\circ-\alpha)\big)^2+ \\[0.2cm]
    &\qquad\;\; +\big(A-B\cos(144^\circ-\alpha)\big)^2+ \\[0.2cm]
    &\qquad +\big(A-B\cos(144^\circ+\alpha)\big)^2+\big(A-B\cos(72^\circ+\alpha)\big)^2= \\[0.2cm]
    & =5A^2-2AB\Big(\cos\alpha+\cos(72^\circ-\alpha)+\cos(144^\circ-\alpha)+ \\[0.2cm]
    &\qquad\qquad\qquad\qquad +\cos(144^\circ+\alpha)+\cos(72^\circ+\alpha)\Big)+ \\[0.2cm]
    &\qquad\quad +B^2\Big(\cos^2\alpha+\cos^2(72^\circ-\alpha)+\cos^2(144^\circ-\alpha)+ \\[0.2cm]
    &\qquad\quad +\cos^2(144^\circ+\alpha)+\cos^2(72^\circ+\alpha)\Big).
\end{align*}
\begin{gather*}
    \cos\alpha+\cos(72^\circ-\alpha)+\cos(144^\circ-\alpha)+\cos(144^\circ+\alpha)+  \\[0.2cm]
    +\cos(72^\circ+\alpha)=\cos\alpha+2\cos\alpha\cdot 72^\circ+2\cos\alpha\cdot\cos 144^\circ= \\[0.2cm]
    =\cos\alpha\big(1+2(\cos 72^\circ+\cos 144^\circ)\big).
\end{gather*}
 }
\vskip+0.2cm
\noindent
Indeed,
\vskip+0.02cm
{\fontsize{14}{16pt}\selectfont
\begin{align*}
    \sin 18^\circ & =\frac{\sqrt{5}-1}{4}\,, \\[0.5cm]
    \cos 72^\circ+\cos 144^\circ & =\sin 18^\circ-\cos 36^\circ= \\[0.2cm]
    & =\sin 18^\circ-1(1-2\sin^2 18^\circ)= \\[0.2cm]
    & =\frac{\sqrt{5}-1}{4}-1+\frac{(\sqrt{5}-1)^2}{2}=-\frac{1}{2}\,,
\end{align*}
 }
\vskip+0.2cm
\noindent
i.e.
\vskip+0.02cm
{\fontsize{14}{16pt}\selectfont
\begin{multline*}
    \cos\alpha+\cos(72^\circ-\alpha)+\cos(144^\circ-\alpha)+ \\[0.2cm]
    +\cos(144^\circ+\alpha)+\cos(72^\circ+\alpha)=0.
\end{multline*}
 }
\vskip+0.2cm
\noindent
Calculate the sum of the second powers:
\vskip+0.02cm
{\fontsize{14}{16pt}\selectfont
\allowdisplaybreaks
\begin{align*}
    \cos^2\alpha & +\cos^2(72^\circ-\alpha)+\cos^2(144^\circ-\alpha)+\qquad\qquad\qquad\qquad \\[0.2cm]
    &\qquad\qquad +\cos^2(144^\circ+\alpha)+\cos^2(72^\circ+\alpha)= \\[0.2cm]
    & =\frac{1}{2}\,\Big(5+\cos 2\alpha+\cos(144^\circ-2\alpha)+\cos(144^\circ+2\alpha)+ \\[0.2cm]
    &\qquad\qquad +\cos(288^\circ-2\alpha)+\cos(288^\circ+2\alpha)\Big)= \\[0.2cm]
    & =\frac{1}{2}\,\Big(5+\cos 2\alpha+\cos(144^\circ-2\alpha)+\cos(144^\circ+2\alpha)+ \\[0.2cm]
    &\qquad\qquad +\cos(72^\circ+2\alpha)+\cos(72^\circ-2\alpha)\Big)=\frac{5}{2}\,.
\end{align*}
 }
\vskip+0.2cm
\noindent
So,
\vskip+0.02cm
{\fontsize{14}{16pt}\selectfont
\begin{multline*}
    S_5^{(4)}=5A^2+\frac{5}{2}\,B^2= \\[0.2cm]
    =5(r^2+\ell^2)^2+10r^2\ell^2=5(r^4+\ell^4+4r^2\ell^2).
\end{multline*}
 }

\bigskip
\begin{theorem}\label{th:5}
{\fontsize{14}{16pt}\selectfont
The locus of points, such that the sum of the fourth powers {\fontsize{14}{16pt}\selectfont $S_5^{(4)}$} of the distances to the vertexes of fixed regular pentagon is constant, is a circle if
\vskip+0.02cm
{\fontsize{14}{16pt}\selectfont
$$  S_5^{(4)}>5r^4,     $$
 }
\vskip+0.2cm
\noindent
where {\fontsize{14}{16pt}\selectfont $r$} is the distance from the center of the regular pentagon to the vertex. The center of the circle is at the center of the regular pentagon and the radius {\fontsize{14}{16pt}\selectfont $\ell$} satisfies the following condition:
\vskip+0.02cm
{\fontsize{14}{16pt}\selectfont
$$  S_5^{(4)}=5\Big[(r^2+\ell^2)^2+2(r\ell)^2\Big].     $$
 }
\vskip+0.2cm

If {\fontsize{14}{16pt}\selectfont $S_5^{(4)}=5r^4$} the locus is a point -- the center of the pentagon and if {\fontsize{14}{16pt}\selectfont $S_5^{(4)}<5r^4$} it is the empty set. }
\end{theorem}

\pagebreak
{\fontsize{14}{16pt}\selectfont
\allowdisplaybreaks
\begin{align*}
    S_5^{(6)} & =(A-B\cos\alpha)^3+\big(A-B\cos(72^\circ-\alpha)\big)^3+ \\[0.2cm]
    &\qquad +\big(A-B\cos(144^\circ-\alpha)\big)^3+\big(A-B\cos(144^\circ+\alpha)\big)^3+ \\[0.2cm]
    &\qquad +\big(A-B\cos(72^\circ+\alpha)\big)^3= \\[0.2cm]
    & =A^3-3A^2B\cos\alpha+3AB^2\cos^2\alpha-B^3\cos^3\alpha+ \\[0.2cm]
    &\qquad +A^3-3A^2B\cos(72^\circ-\alpha)+3AB^2\cos^2(72^\circ-\alpha)- \\[0.2cm]
    &\qquad\qquad\qquad\qquad -B^3\cos^3(72^\circ-\alpha)+ \\[0.2cm]
    &\qquad +A^3-3A^2B\cos(144^\circ-\alpha)+3AB^2\cos^2(144^\circ-\alpha)- \\[0.2cm]
    &\qquad\qquad\qquad -B^3\cos^3(144^\circ-\alpha)+ \\[0.2cm]
    &\qquad +A^3-3A^2B\cos(144^\circ+\alpha)+3AB^2\cos^2(144^\circ+\alpha)- \\[0.2cm]
    &\qquad\qquad\qquad\qquad -B^3\cos^3(144^\circ+\alpha)+ \\[0.2cm]
    &\qquad +A^3-3A^2B\cos(72^\circ+\alpha)+3AB^2\cos^2(72^\circ+\alpha)- \\[0.2cm]
    &\qquad\qquad\qquad\qquad -B^3\cos^3(72^\circ+\alpha)= \\[0.2cm]
    & =5A^3-3A^2B\Big(\cos\alpha+\cos(72^\circ-\alpha)+\cos(144^\circ-\alpha)+ \\[0.2cm]
    &\qquad\qquad\qquad\qquad +\cos(144^\circ+\alpha)+\cos(72^\circ+\alpha)\Big)+ \\[0.2cm]
    &\qquad +3AB^2\Big(\cos^2\alpha+\cos^2(72^\circ-\alpha)+\cos^2(144^\circ-\alpha)+ \\[0.2cm]
    &\qquad\qquad\qquad\qquad +\cos^2(144^\circ+\alpha)+\cos^2(72^\circ+\alpha)\Big)- \\[0.2cm]
    &\qquad -B^3\Big(\cos^3\alpha+\cos^3(72^\circ-\alpha)+\cos^3(144^\circ-\alpha)+ \\[0.2cm]
    &\qquad\qquad\qquad\qquad +\cos^3(144^\circ+\alpha)+\cos^3(72^\circ+\alpha)\Big).
\end{align*}
 }
\vskip+0.2cm
\noindent
From the preceding part we know the values of the sums of the first and the second powers. Calculate the sum of the third powers:
\vskip+0.02cm
{\fontsize{14}{16pt}\selectfont
\allowdisplaybreaks
\begin{align*}
    \cos^3\alpha & +\cos^3(72^\circ-\alpha)+\cos^3(144^\circ-\alpha)+ \qquad\qquad\qquad \\[0.2cm]
    &\qquad\qquad +\cos^3(144^\circ+\alpha)+\cos^3(72^\circ+\alpha)= \\[0.2cm]
    & =\frac{1}{4}\,\Big(3\cos\alpha+\cos 3\alpha+ \\[0.2cm]
    &\qquad\qquad +3\cos(72^\circ-\alpha)+\cos(216^\circ-3\alpha)+ \\[0.2cm]
    &\qquad\qquad +3\cos(144^\circ-\alpha)+\cos(432^\circ-3\alpha)+ \\[0.2cm]
    &\qquad\qquad +3\cos(144^\circ+\alpha)+\cos(432^\circ+3\alpha)+ \\[0.2cm]
    &\qquad\qquad +3\cos(72^\circ+\alpha)+\cos(216^\circ+3\alpha)\Big)= \\[0.2cm]
    & =\frac{1}{4}\,\Big(3\Big(\cos\alpha+\cos(72^\circ-\alpha)+\cos(144^\circ-\alpha)+  \\[0.2cm]
    &\qquad +\cos(144^\circ+\alpha)+\cos(72^\circ+\alpha)\Big)+ \\[0.2cm]
    &\qquad +\cos 3\alpha+\cos(216^\circ-3\alpha)+\cos(432^\circ-3\alpha)+ \\[0.2cm]
    &\qquad +3\cos(432^\circ+3\alpha)+\cos(216^\circ+3\alpha)\Big)= \\[0.2cm]
    & =\frac{1}{4}\,\Big(\cos 3\alpha+\cos(144^\circ+3\alpha)+\cos(72^\circ-3\alpha)+ \\[0.2cm]
    &\qquad +\cos(72^\circ+3\alpha)+\cos(144^\circ-3\alpha)\Big)=0.
\end{align*}
 }
\vskip+0.2cm
\noindent
So,
\vskip+0.02cm
{\fontsize{14}{16pt}\selectfont
\begin{multline*}
    S_5^{(6)}=5A^3+3AB^2\cdot\frac{5}{2}= \\[0.2cm]
    =5A\Big(A^2+\frac{3B^2}{2}\Big)=5(r^2+\ell^2)(r^4+\ell^4+8r^2\ell^2).
\end{multline*}
 }

\bigskip
\begin{theorem}\label{th:6}
{\fontsize{14}{16pt}\selectfont
The locus of points, such that the sum of the sixth powers {\fontsize{14}{16pt}\selectfont $S_5^{(6)}$} of the distances to the vertexes of fixed regular pentagon is constant, is a circle if
\vskip+0.02cm
{\fontsize{14}{16pt}\selectfont
$$  S_5^{(6)}>5r^6,     $$
 }
\vskip+0.2cm
\noindent
where {\fontsize{14}{16pt}\selectfont $r$} is the distance from the center of the regular pentagon to the vertex. The center of the circle is at the center of the regular pentagon and the radius {\fontsize{14}{16pt}\selectfont $\ell$} satisfies the following condition:
\vskip+0.02cm
{\fontsize{14}{16pt}\selectfont
$$  S_5^{(6)}=5\Big[(r^2+\ell^2)^3+6(r^2+\ell^2)(r\ell)^2\Big].     $$
 }
\vskip+0.2cm

If {\fontsize{14}{16pt}\selectfont $S_5^{(6)}=5r^6$} the locus is a point -- the center of the regular pentagon and if {\fontsize{14}{16pt}\selectfont $S_5^{(6)}<5r^6$} it is the empty set. }
\end{theorem}

\vskip+0.02cm
{\fontsize{14}{16pt}\selectfont
\allowdisplaybreaks
\begin{align*}
    S_5^{(8)} & =(A-B\cos\alpha)^4+\big(A-B\cos(72^\circ-\alpha)\big)^4+ \\[0.2cm]
    &\quad +\big(A-B\cos(144^\circ-\alpha)\big)^4+\big(A-B\cos(144^\circ+\alpha)\big)^4+ \\[0.2cm]
    &\quad +\big(A-B\cos(72^\circ+\alpha)\big)^4= 
\end{align*}
\pagebreak
\begin{align*}
    & =A^4-4A^3B\cos\alpha+6A^2B^2\cos^2\alpha-  \\[0.2cm]
    &\qquad\qquad\qquad -4AB^3\cos^3\alpha+B^4\cos^4\alpha+\cdots= \\[0.2cm]
    & =5A^4-4A^3B\Big(\cos\alpha+\cos(72^\circ-\alpha)+\cos(144^\circ-\alpha)+ \\[0.2cm]
    &\qquad\qquad +\cos(144^\circ+\alpha)+\cos(72^\circ+\alpha)\Big)+ \\[0.2cm]
    &\quad +6A^2B^2\Big(\cos^2\alpha+\cos^2(72^\circ-\alpha)+\cos^2(144^\circ-\alpha)+ \\[0.2cm]
    &\qquad\qquad\qquad +\cos^2(144^\circ+\alpha)+\cos^2(72^\circ+\alpha)\Big)- \\[0.2cm]
    &\quad -4AB^3\Big(\cos^3\alpha+\cos^3(72^\circ-\alpha)+\cos^3(144^\circ-\alpha)+ \\[0.2cm]
    &\qquad\qquad\qquad +\cos^3(144^\circ+\alpha)+\cos^3(72^\circ+\alpha)\Big)+ \\[0.2cm]
    &\quad+B^4\Big(\cos^4\alpha+\cos^4(72^\circ-\alpha)+\cos^4(144^\circ-\alpha)+ \\[0.2cm]
    &\qquad\qquad\qquad +\cos^4(144^\circ+\alpha)+\cos^4(72^\circ+\alpha)\Big).
\end{align*}
 }
\vskip+0.2cm

Calculate the sum of the fourth powers:
\vskip+0.02cm
{\fontsize{14}{16pt}\selectfont
\allowdisplaybreaks
\begin{align*}
    \cos^4\alpha & +\cos^4(72^\circ-\alpha)+\cos^4(144^\circ-\alpha)+ \\[0.3cm]
    &\qquad\qquad +\cos^4(144^\circ+\alpha)+\cos^4(72^\circ+\alpha)= \\[0.3cm]
    & =\frac{1}{8}\,\Big(\cos4\alpha+8\cos^2\alpha-1+  \\[0.3cm]
    &\qquad +\cos(288^\circ-4\alpha)+8\cos^2(72^\circ-\alpha)-1+  \\[0.3cm]
    &\qquad +\cos(576^\circ-4\alpha)+8\cos^2(144^\circ-\alpha)-1+  \\[0.3cm]
    &\qquad +\cos(576^\circ+4\alpha)+8\cos^2(144^\circ+\alpha)-1+  \\[0.3cm]
    &\qquad +\cos(288^\circ+4\alpha)+8\cos^2(72^\circ+\alpha)-1\Big)= \\[0.3cm]
    & =\frac{1}{8}\,\Big(-5+8\Big(\cos^2\alpha+\cos^2(72^\circ-\alpha)+ \\[0.3cm]
    &\quad +\cos^2(144^\circ-\alpha)+\cos^2(144^\circ+\alpha)+\cos^2(72^\circ+\alpha)\Big)+ \\[0.3cm]
    &\qquad +\cos4\alpha+\cos(288^\circ-4\alpha)+\cos(576^\circ-4\alpha)+ \\[0.3cm]
    &\qquad\qquad +\cos(576^\circ+4\alpha)+\cos(288^\circ+4\alpha)\Big)= \\[0.3cm]
    & =\frac{1}{8}\,\Big(-5+8\cdot\frac{5}{2}+\cos 4\alpha+\cos(72^\circ+4\alpha)+ \\[0.3cm]
    &\quad +\cos(144^\circ+4\alpha)+\cos(144^\circ-4\alpha)+\cos(72^\circ-4\alpha)\Big)= \\
    & =\frac{15}{8}\,, \\[0.5cm]
    S_5^{(8)} & =5A^4+6A^2B^2\cdot\frac{5}{2}+B^4\cdot\frac{15}{8}= \\[0.3cm]
    & =5(r^2+\ell^2)^2(r^4+\ell^4+14r^2\ell^2)+30r^4\ell^4.
\end{align*}
 }

\bigskip
\begin{theorem}\label{th:7}
{\fontsize{14}{16pt}\selectfont
The locus of points, such that the sum of the eighth powers {\fontsize{14}{16pt}\selectfont $S_5^{(8)}$} of the distances to the vertexes of fixed regular pentagon is constant, is a circle if
\vskip+0.02cm
{\fontsize{14}{16pt}\selectfont
$$  S_5^{(8)}>5r^8,     $$
 }
\vskip+0.2cm
\noindent
where {\fontsize{14}{16pt}\selectfont $r$} is the distance from the center of the regular pentagon to the vertex. The center of the circle is at the center of the regular pentagon and the radius {\fontsize{14}{16pt}\selectfont $\ell$} satisfies the following condition:
\vskip+0.02cm
{\fontsize{14}{16pt}\selectfont
$$  S_5^{(8)}=5\Big[(r^2+\ell^2)^4+12(r^2+\ell^2)^2(r\ell)^2+6(r\ell)^4\Big].     $$
 }
\vskip+0.2cm

If {\fontsize{14}{16pt}\selectfont $S_5^{(8)}=5r^8$} the locus is a point -- the center of the regular pentagon and if {\fontsize{14}{16pt}\selectfont $S_5^{(8)}<5r^8$} it is the empty set. }
\end{theorem}

\bigskip
\noindent \textbf{Note.} The sum of the tenth powers contains {\fontsize{14}{16pt}\selectfont $\alpha$}, so the locus is not a circle.

\vskip+1.5cm
\section{\fontsize{18}{20pt}\selectfont \textbf{Circle as Locus of Constant Sum of Powers}}
\vskip+1cm

Summarize obtained results. If {\fontsize{14}{16pt}\selectfont $d(X,P_i)$} is the distance from {\fontsize{14}{16pt}\selectfont $X$} to the vertex {\fontsize{14}{16pt}\selectfont $P_i$} of the regular {\fontsize{14}{16pt}\selectfont $n$}-sided polygon and {\fontsize{14}{16pt}\selectfont $r$} is the distance from the center of the regular polygon to the vertex, we have:

\bigskip
\bigskip
\textbf{1.} {\fontsize{14}{16pt}\selectfont $\displaystyle \Big\{X:\;\;\sum_{i=1}^n d^2(X,P_i)=S_n^{(2)}=const\Big\}=$}
\vskip+0.2cm
{\fontsize{14}{16pt}\selectfont
$$  \qquad\qquad\qquad =\begin{cases}
                    \text{\fontsize{12}{14pt}\selectfont circle}, & \text{\fontsize{12}{14pt}\selectfont if}\;\; S_n^{(2)}>nr^2 \\[0.15cm]
                    \text{\fontsize{12}{14pt}\selectfont point}, & \text{\fontsize{12}{14pt}\selectfont if}\;\; S_n^{(2)}=nr^2 \\[0.15cm]
                    \;\;\varnothing\;\;, & \text{\fontsize{12}{14pt}\selectfont if}\;\; S_n^{(2)}<nr^2
                \end{cases}\,,     $$
}

\pagebreak
\bigskip
\bigskip
\textbf{2.} {\fontsize{14}{16pt}\selectfont $\displaystyle \Big\{X:\;\;\sum_{i=1}^3 d^4(X,P_i)=S_3^{(4)}=const\Big\}=$}
\vskip+0.2cm
{\fontsize{14}{16pt}\selectfont
$$  \qquad\qquad\qquad =\begin{cases}
                    \text{\fontsize{12}{14pt}\selectfont circle}, & \text{\fontsize{12}{14pt}\selectfont if}\;\; S_3^{(4)}>3r^4 \\[0.15cm]
                    \text{\fontsize{12}{14pt}\selectfont point}, & \text{\fontsize{12}{14pt}\selectfont if}\;\; S_3^{(4)}=3r^4 \\[0.15cm]
                    \;\;\varnothing\;\;, & \text{\fontsize{12}{14pt}\selectfont if}\;\; S_3^{(4)}<3r^4
                \end{cases}\,.     $$
}

\bigskip
\bigskip
\textbf{3.} {\fontsize{14}{16pt}\selectfont $\displaystyle \Big\{X:\;\;\sum_{i=1}^4 d^4(X,P_i)=S_4^{(4)}=const\Big\}=$}
\vskip+0.2cm
{\fontsize{14}{16pt}\selectfont
$$  \qquad\qquad\qquad =\begin{cases}
                    \text{\fontsize{12}{14pt}\selectfont circle}, & \text{\fontsize{12}{14pt}\selectfont if}\;\; S_4^{(4)}>4r^4 \\[0.15cm]
                    \text{\fontsize{12}{14pt}\selectfont point}, & \text{\fontsize{12}{14pt}\selectfont if}\;\; S_4^{(4)}=4r^4 \\[0.15cm]
                    \;\;\varnothing\;\;, & \text{\fontsize{12}{14pt}\selectfont if}\;\; S_4^{(4)}<4r^4
                \end{cases}\,.     $$
}

\bigskip
\bigskip
\textbf{4.} {\fontsize{14}{16pt}\selectfont $\displaystyle \Big\{X:\;\;\sum_{i=1}^4 d^6(X,P_i)=S_4^{(6)}=const\Big\}=$}
\vskip+0.2cm
{\fontsize{14}{16pt}\selectfont
$$  \qquad\qquad\qquad =\begin{cases}
                    \text{\fontsize{12}{14pt}\selectfont circle}, & \text{\fontsize{12}{14pt}\selectfont if}\;\; S_4^{(6)}>4r^6 \\[0.15cm]
                    \text{\fontsize{12}{14pt}\selectfont point}, & \text{\fontsize{12}{14pt}\selectfont if}\;\; S_4^{(6)}=4r^6 \\[0.15cm]
                    \;\;\varnothing\;\;, & \text{\fontsize{12}{14pt}\selectfont if}\;\; S_4^{(6)}<4r^6
                \end{cases}\,.     $$
}

\pagebreak
\bigskip
\bigskip
\textbf{5.} {\fontsize{14}{16pt}\selectfont $\displaystyle \Big\{X:\;\;\sum_{i=1}^5 d^4(X,P_i)=S_5^{(4)}=const\Big\}=$}
\vskip+0.2cm
{\fontsize{14}{16pt}\selectfont
$$  \qquad\qquad\qquad =\begin{cases}
                    \text{\fontsize{12}{14pt}\selectfont circle}, & \text{\fontsize{12}{14pt}\selectfont if}\;\; S_5^{(4)}>5r^4 \\[0.15cm]
                    \text{\fontsize{12}{14pt}\selectfont point}, & \text{\fontsize{12}{14pt}\selectfont if}\;\; S_5^{(4)}=5r^4 \\[0.15cm]
                    \;\;\varnothing\;\;, & \text{\fontsize{12}{14pt}\selectfont if}\;\; S_5^{(4)}<5r^4
                \end{cases}\,.     $$
}

\bigskip
\bigskip
\textbf{6.} {\fontsize{14}{16pt}\selectfont $\displaystyle \Big\{X:\;\;\sum_{i=1}^5 d^6(X,P_i)=S_5^{(6)}=const\Big\}=$}
\vskip+0.2cm
{\fontsize{14}{16pt}\selectfont
$$  \qquad\qquad\qquad =\begin{cases}
                    \text{\fontsize{12}{14pt}\selectfont circle}, & \text{\fontsize{12}{14pt}\selectfont if}\;\; S_5^{(6)}>5r^6 \\[0.15cm]
                    \text{\fontsize{12}{14pt}\selectfont point}, & \text{\fontsize{12}{14pt}\selectfont if}\;\; S_5^{(6)}=5r^6 \\[0.15cm]
                    \;\;\varnothing\;\;, & \text{\fontsize{12}{14pt}\selectfont if}\;\; S_5^{(6)}<5r^6
                \end{cases}\,.    $$
}

\bigskip
\bigskip
\textbf{7.} {\fontsize{14}{16pt}\selectfont $\displaystyle \Big\{X:\;\;\sum_{i=1}^5 d^8(X,P_i)=S_5^{(8)}=const\Big\}=$}
\vskip+0.2cm
{\fontsize{14}{16pt}\selectfont
$$  \qquad\qquad\qquad =\begin{cases}
                    \text{\fontsize{12}{14pt}\selectfont circle}, & \text{\fontsize{12}{14pt}\selectfont if}\;\; S_5^{(8)}>5r^8 \\[0.15cm]
                    \text{\fontsize{12}{14pt}\selectfont point}, & \text{\fontsize{12}{14pt}\selectfont if}\;\; S_5^{(8)}=5r^8 \\[0.15cm]
                    \;\;\varnothing\;\;, & \text{\fontsize{12}{14pt}\selectfont if}\;\; S_5^{(8)}<5r^8
                \end{cases}\,.     $$
}

\bigskip

Now we are ready to state general hypothesis:

\pagebreak
\begin{hypothesis}
{\fontsize{14}{16pt}\selectfont
The locus of points, such that the sum of the   \linebreak  {\fontsize{14}{16pt}\selectfont $(2m)$}-th powers {\fontsize{14}{16pt}\selectfont $S_n^{(2m)}$} of the distances to the vertexes of fixed re\-gu\-lar {\fontsize{14}{16pt}\selectfont $n$}-sided polygon is constant, is a circle if
\vskip+0.02cm
{\fontsize{14}{16pt}\selectfont
$$  S_n^{(2m)}>nr^{2m}, \;\;\text{where}\;\; m=1,2,\dots,n-1        $$
 }
\vskip+0.2cm
\noindent
and {\fontsize{14}{16pt}\selectfont $r$} is the distance from the center of the regular polygon to the vertex. The center of the circle is at the center of the regular polygon.

If {\fontsize{14}{16pt}\selectfont $S_n^{(2m)}=nr^{2m}$} the locus is a point -- the center of the polygon and if {\fontsize{14}{16pt}\selectfont $S_n^{(2m)}<nr^{2m}$} it is the empty set. }
\end{hypothesis}
\bigskip
\bigskip

Thus, for each even {\fontsize{14}{16pt}\selectfont $(2m)$}-th power, so that {\fontsize{14}{16pt}\selectfont $m<n$}:
\vskip+0.02cm
{\fontsize{14}{16pt}\selectfont
\begin{multline*}
    \qquad \Big\{X:\;\;\sum_{i=1}^n d^{2m}(X,P_i)=S_n^{(2m)}=const\Big\}= \\
    =\begin{cases}
                    \text{\fontsize{12}{14pt}\selectfont circle}, & \text{\fontsize{12}{14pt}\selectfont if}\;\; S_n^{(2m)}>nr^{2m} \\[0.15cm]
                    \text{\fontsize{12}{14pt}\selectfont point}, & \text{\fontsize{12}{14pt}\selectfont if}\;\; S_n^{(2m)}=nr^{2m} \\[0.15cm]
                    \;\;\varnothing\;\;, & \text{\fontsize{12}{14pt}\selectfont if}\;\; S_n^{(2m)}<nr^{2m}
                \end{cases}\,.  \qquad
\end{multline*}
}
\vskip+0.2cm
\noindent

To prove the hypothesis, first of all two lemmas must be proved.

\newpage
\section{\fontsize{18}{20pt}\selectfont \textbf{Cosine Multiple Arguments \\ and Cosine Powers Sums}}
\vskip+1cm

\begin{lemma}\label{lem:1}
{\fontsize{14}{16pt}\selectfont
For arbitrary positive integers {\fontsize{14}{16pt}\selectfont $m$} and {\fontsize{14}{16pt}\selectfont $n$}, so that   \linebreak   {\fontsize{14}{16pt}\selectfont $m<n$}, the following condition
\vskip+0.02cm
{\fontsize{14}{16pt}\selectfont
$$  \sum_{k=1}^n \cos\bigg(m\Big(\alpha-(k-1)\,\frac{360^{\circ}}{n}\Big)\bigg)=0        $$
 }
\vskip+0.2cm
\noindent
is satisfied, where {\fontsize{14}{16pt}\selectfont $\alpha$} is an arbitrary angle. }
\end{lemma}

Denote by
\vskip+0.02cm
{\fontsize{14}{16pt}\selectfont
$$  P\!=\!e^{\text{\fontsize{12}{14pt}\selectfont $im\alpha$}}+e^{\text{\fontsize{12}{14pt}\selectfont $im(\alpha-\frac{360^{\circ}}{n})$}}+
        e^{\text{\fontsize{12}{14pt}\selectfont $im(\alpha-2\,\frac{360^{\circ}}{n})$}}+\cdots+
                e^{\text{\fontsize{12}{14pt}\selectfont $im(\alpha-(n-1)\,\frac{360^{\circ}}{n})$}}.     $$
 }
\vskip+0.2cm
\noindent
By using Euler's formula
\vskip+0.02cm
{\fontsize{14}{16pt}\selectfont
$$  e^{i\varphi}=\cos\varphi+i\sin\varphi      $$
 }
\vskip+0.2cm
\noindent
real part of {\fontsize{14}{16pt}\selectfont $P$} is:
\vskip+0.02cm
{\fontsize{14}{16pt}\selectfont
$$  \RRe(P)=\sum_{k=1}^n \cos\bigg(m\Big(\alpha-(k-1)\,\frac{360^{\circ}}{n}\Big)\bigg).      $$
 }
\vskip+0.2cm

The formula of the sum of geometric progression gives:
\vskip+0.02cm
{\fontsize{14}{16pt}\selectfont
\begin{gather*}
\begin{aligned}
    P & =e^{\text{\fontsize{12}{14pt}\selectfont $im\alpha$}}
        \bigg(1\!+\!e^{\text{\fontsize{12}{14pt}\selectfont $-im\,\frac{360^{\circ}}{n}$}}\!+\!
            \big(e^{\text{\fontsize{12}{14pt}\selectfont $-im\,\frac{360^{\circ}}{n}$}}\big)^2\!+\cdots+\!
                    \big(e^{\text{\fontsize{12}{14pt}\selectfont $-im\,\frac{360^{\circ}}{n}$}}\big)^{n-1}\bigg)= \\[0.2cm]
    & =e^{\text{\fontsize{12}{14pt}\selectfont $im\alpha$}}\,
            \frac{1-(e^{\text{\fontsize{12}{14pt}\selectfont $-im\,\frac{360^{\circ}}{n}$}})^n}
                        {1-e^{\text{\fontsize{12}{14pt}\selectfont $-im\,\frac{360^{\circ}}{n}$}}}\,,
\end{aligned} \\[0.5cm]
    e^{-im\,360^{\circ}}=\cos(-360^{\circ}m)+i(-360^{\circ}m)=1.
\end{gather*}
 }
\vskip+0.2cm

Since {\fontsize{14}{16pt}\selectfont $m<n$},\; {\fontsize{14}{16pt}\selectfont $e^{-im\,\frac{360^{\circ}}{n}}\neq 1$}. So {\fontsize{14}{16pt}\selectfont $P=0$}, i.e. {\fontsize{14}{16pt}\selectfont $\RRe(P)=0$} which proves the Lemma \ref{lem:1}.

\bigskip
\begin{note}
{\fontsize{14}{16pt}\selectfont If {\fontsize{14}{16pt}\selectfont $m=n$}, the sum equals\; {\fontsize{14}{16pt}\selectfont $n\cos(n\alpha)$}. The sum is not zero only in cases when {\fontsize{14}{16pt}\selectfont $m$} is multiple of {\fontsize{14}{16pt}\selectfont $n$}. If {\fontsize{14}{16pt}\selectfont $m=pn$}, the sum is\;  {\fontsize{14}{16pt}\selectfont $n\cos(pn\alpha)$}. }
\end{note}

\bigskip

\begin{lemma}\label{lem:2}
{\fontsize{14}{16pt}\selectfont For arbitrary positive integers {\fontsize{14}{16pt}\selectfont $m$} and {\fontsize{14}{16pt}\selectfont $n$}, so that   \linebreak  {\fontsize{14}{16pt}\selectfont $m<n$} and for an arbitrary angle {\fontsize{14}{16pt}\selectfont $\alpha$} the following conditions are satisfied:
\begin{enumerate}
\item[] if {\fontsize{14}{16pt}\selectfont $m$} is odd
\vskip+0.02cm
{\fontsize{14}{16pt}\selectfont
$$  \sum_{k=1}^n \cos^m\Big(\alpha-(k-1)\,\frac{360^{\circ}}{n}\Big)=0;     $$
 }
\vskip+0.2cm

\item[] if {\fontsize{14}{16pt}\selectfont $m$} is even
{\fontsize{14}{16pt}\selectfont
$$  \sum_{k=1}^n \cos^m\Big(\alpha-(k-1)\,\frac{360^{\circ}}{n}\Big)=n\,
                \frac{\Big(\!\!\begin{array}{c} \text{\fontsize{13}{14pt}\selectfont $m$} \\ \text{\fontsize{13}{14pt}\selectfont $\frac{m}{2}$} \end{array}\!\!\Big)}
                            {2^m}\,.     $$
 }
\end{enumerate}
}
\end{lemma}
\bigskip
\bigskip

By using power-reduction formula for cosine, when {\fontsize{14}{16pt}\selectfont $m$} is odd:
\vskip+0.02cm
{\fontsize{14}{16pt}\selectfont
$$  \cos^m\theta=\frac{2}{2^m} \sum_{k=0}^{\frac{m-1}{2}}
        \left(\!\!\!\begin{array}{c} \text{\fontsize{13}{14pt}\selectfont $m$} \\ \text{\fontsize{13}{14pt}\selectfont $k$} \end{array}\!\!\!\right)
                            \cos\big((m-2k)\theta\big),        $$
 }
\vskip+0.2cm
\noindent
we have
\vskip+0.02cm
{\fontsize{14}{16pt}\selectfont
\allowdisplaybreaks
\begin{align*}
    \sum_{k=1}^n & \cos^m\Big(\alpha-(k-1)\,\frac{360^{\circ}}{n}\Big)= \\[0.2cm]
    & =\cos^m\alpha+\cos^m\Big(\alpha-\frac{360^{\circ}}{n}\Big)+\cdots+ \\[0.2cm]
    &\qquad\qquad\qquad\qquad\qquad +\cos^m\Big(\alpha-(n-1)\,\frac{360^{\circ}}{n}\Big)= 
\end{align*}
\begin{align*}
    & =\frac{2}{2^m}\,\Bigg[\left(\!\!\!\begin{array}{c} \text{\fontsize{13}{14pt}\selectfont $m$} \\ \text{\fontsize{13}{14pt}\selectfont $0$} \end{array}\!\!\!\right)
                                    \cos m\alpha+
        \left(\!\!\!\begin{array}{c} \text{\fontsize{13}{14pt}\selectfont $m$} \\ \text{\fontsize{13}{14pt}\selectfont $1$} \end{array}\!\!\!\right)
                \cos (m-2)\alpha+\cdots+
            \left(\!\!\!\begin{array}{c} \text{\fontsize{13}{14pt}\selectfont $m$} \\ \text{\fontsize{13}{14pt}\selectfont $\frac{m-1}{2}$} \end{array}\!\!\!\right)
                                    \cos\alpha+ \\[0.4cm]
    &\quad +\left(\!\!\!\begin{array}{c} \text{\fontsize{13}{14pt}\selectfont $m$} \\ \text{\fontsize{13}{14pt}\selectfont $0$} \end{array}\!\!\!\right)
                    \cos m\Big(\alpha-\frac{360^{\circ}}{n}\Big)+
            \left(\!\!\!\begin{array}{c} \text{\fontsize{13}{14pt}\selectfont $m$} \\ \text{\fontsize{13}{14pt}\selectfont $1$} \end{array}\!\!\!\right)
                        \cos(m-2)\Big(\alpha-\frac{360^{\circ}}{n}\Big)+\cdots+ \\[0.4cm]
    &\qquad\qquad\qquad\qquad +\left(\!\!\!\begin{array}{c} \text{\fontsize{13}{14pt}\selectfont $m$} \\ \text{\fontsize{13}{14pt}\selectfont $\frac{m-1}{2}$} \end{array}\!\!\!\right)
                    \cos\Big(\alpha-\frac{360^{\circ}}{n}\Big)+\;\cdots\;+ \\[0.4cm]
    &\quad +\left(\!\!\!\begin{array}{c} \text{\fontsize{13}{14pt}\selectfont $m$} \\ \text{\fontsize{13}{14pt}\selectfont $0$} \end{array}\!\!\!\right)
                        \cos m\Big(\alpha\!-\!(n\!-\!1)\,\frac{360^{\circ}}{n}\Big)\!+\!
            \left(\!\!\!\begin{array}{c} \text{\fontsize{13}{14pt}\selectfont $m$} \\ \text{\fontsize{13}{14pt}\selectfont $1$} \end{array}\!\!\!\right)
                    \cos(m\!-\!2)\Big(\alpha\!-\!(n\!-\!1)\,\frac{360^{\circ}}{n}\Big)\!+ \\[0.4cm]
    &\qquad\qquad\quad\;\;\; +\cdots+\left(\!\!\!\begin{array}{c} \text{\fontsize{13}{14pt}\selectfont $m$} \\ \text{\fontsize{13}{14pt}\selectfont $\frac{m-1}{2}$} \end{array}\!\!\!\right)
                        \cos\Big(\alpha-(n-1)\,\frac{360^{\circ}}{n}\Big)\Bigg]= 
\end{align*}
\pagebreak
\begin{align*}
    & =\frac{2}{2^m}\,\Bigg[\left(\!\!\!\begin{array}{c} \text{\fontsize{13}{14pt}\selectfont $m$} \\ \text{\fontsize{13}{14pt}\selectfont $0$} \end{array}\!\!\!\right)
                \bigg(\cos m\alpha+\cos m\Big(\alpha-\frac{360^{\circ}}{n}\Big)+\cdots+ \\[0.4cm]
    &\qquad\qquad\qquad\qquad\qquad +\cos m\Big(\alpha-(n-1)\,\frac{360^{\circ}}{n}\Big)\bigg)+ \\[0.4cm]
    &\qquad +\left(\!\!\!\begin{array}{c} \text{\fontsize{13}{14pt}\selectfont $m$} \\ \text{\fontsize{13}{14pt}\selectfont $1$} \end{array}\!\!\!\right)
                \bigg(\cos(m\!-\!2)\alpha\!+\!\cos(m\!-\!2)\Big(\alpha\!-\!\frac{360^{\circ}}{n}\Big)\!+\cdots+ \\[0.4cm]
    &\qquad\qquad\qquad\qquad\qquad +\cos(m-2)\Big(\alpha-(n-1)\,\frac{360^{\circ}}{n}\Big)\bigg)+ \\[0.4cm]
    &\qquad\qquad \qquad\qquad \vdots \\[0.4cm]
    &\qquad +\left(\!\!\!\begin{array}{c} \text{\fontsize{13}{14pt}\selectfont $m$} \\ \text{\fontsize{13}{14pt}\selectfont $\frac{m-1}{2}$} \end{array}\!\!\!\right)
                \bigg(\cos\alpha+\cos\Big(\alpha-\frac{360^{\circ}}{n}\Big)+\cdots+ \\[0.4cm]
    &\qquad\qquad\qquad\qquad\qquad +\cos\Big(\alpha-(n-1)\,\frac{360^{\circ}}{n}\Big)\bigg)\Bigg].
\end{align*}
 }
\vskip+0.2cm
\noindent
Since {\fontsize{14}{16pt}\selectfont $m<n$}, from the Lemma~\ref{lem:1} follows -- each sum equals zero, which proves the first part of the Lemma \ref{lem:2}.

When {\fontsize{14}{16pt}\selectfont $m$} is even power-reduction formula for cosine is:
\vskip+0.02cm
{\fontsize{14}{16pt}\selectfont
$$  \cos^m\theta=\frac{1}{2^m}
            \left(\!\!\!\begin{array}{c} \text{\fontsize{13}{14pt}\selectfont $m$} \\ \text{\fontsize{13}{14pt}\selectfont $\frac{m}{2}$} \end{array}\!\!\!\right)+
                        \frac{2}{2^m} \sum_{k=0}^{\frac{m}{2}-1}
            \left(\!\!\!\begin{array}{c} \text{\fontsize{13}{14pt}\selectfont $m$} \\ \text{\fontsize{13}{14pt}\selectfont $k$} \end{array}\!\!\!\right)
                        \cos\big((m-2k)\theta\big).       $$
 }
\vskip+0.2cm
\noindent
By using the same manner, as in preceding part the sum of the second addends gives zero, and because the number of the first addends is {\fontsize{14}{16pt}\selectfont $n$} total sum equals:
\vskip+0.02cm
{\fontsize{14}{16pt}\selectfont
$$  n\,
                \frac{\Big(\!\!\begin{array}{c} \text{\fontsize{13}{14pt}\selectfont $m$} \\ \text{\fontsize{13}{14pt}\selectfont $\frac{m}{2}$} \end{array}\!\!\Big)}
                            {2^m}\,,     $$
 }
\vskip+0.2cm
\noindent
which proves the Lemma \ref{lem:2}.

\vskip+1.5cm
\section{\fontsize{18}{20pt}\selectfont \textbf{General Theorem}}
\vskip+0.7cm

\begin{theorem}\label{th:8.1}
{\fontsize{14}{16pt}\selectfont The locus of points, such that the sum of the {\fontsize{14}{16pt}\selectfont $(2m)$}-th powers {\fontsize{14}{16pt}\selectfont $S_n^{(2m)}$} of distances to the vertexes of fixed regular {\fontsize{14}{16pt}\selectfont $n$}-sided polygon is constant, is a circle if
\vskip+0.02cm
{\fontsize{14}{16pt}\selectfont
$$  S_n^{(2m)}>nr^{2m}, \;\;\text{\fontsize{12}{14pt}\selectfont where}\;\; m=1,2,\dots,n-1     $$
 }
\vskip+0.2cm
\noindent
and {\fontsize{14}{16pt}\selectfont $r$} is the distance from the center of the regular polygon to the vertex. The center of the circle is at the center of the regular polygon and the radius {\fontsize{14}{16pt}\selectfont $\ell$} satisfies the following condition:
\vskip+0.02cm
{\fontsize{14}{16pt}\selectfont
$$  S_n^{(2m)}=n\Bigg[(r^2+\ell^2)^m+\sum_{k=1}^{[\frac{m}{2}]}
            \left(\!\!\!\begin{array}{c} \text{\fontsize{13}{14pt}\selectfont $m$} \\ \text{\fontsize{13}{14pt}\selectfont $2k$} \end{array}\!\!\!\right)
                        (r^2+\ell^2)^{m-2k}(r\ell)^{2k}
            \left(\!\!\!\begin{array}{c} \text{\fontsize{13}{14pt}\selectfont $2k$} \\ \text{\fontsize{13}{14pt}\selectfont $k$} \end{array}\!\!\!\right)\Bigg].        $$
 }
\vskip+0.2cm

If {\fontsize{14}{16pt}\selectfont $S_n^{(2m)}=nr^{2m}$} the locus is a point -- the center of the polygon and if {\fontsize{14}{16pt}\selectfont $S_n^{(2m)}<nr^{2m}$} the locus is the empty set. }
\end{theorem}

\bigskip
\bigskip

For {\fontsize{14}{16pt}\selectfont $(2m)$}-th power sum of the distances from an arbitrary point to the vertexes {\fontsize{14}{16pt}\selectfont $P_i$} of regular {\fontsize{14}{16pt}\selectfont $n$}-sided polygon, we have:
\vskip+0.02cm
{\fontsize{14}{16pt}\selectfont
\allowdisplaybreaks
\begin{align*}
    S_n^{(2m)} & =(A-B\cos\alpha)^m+\bigg(A-B\cos\Big(\frac{360^{\circ}}{n}-\alpha\Big)\bigg)^m+ \\[0.2cm]
    &\qquad +\bigg(A-B\cos\Big(2\cdot\frac{360^{\circ}}{n}-\alpha\Big)\bigg)^m+\cdots+ \\[0.2cm]
    &\qquad\qquad\qquad +\bigg(A-B\cos\Big((n-1)\,\frac{360^{\circ}}{n}-\alpha\Big)\bigg)^m= \\[0.2cm]
    & =nA^m-\left(\!\!\!\begin{array}{c} \text{\fontsize{13}{14pt}\selectfont $m$} \\ \text{\fontsize{13}{14pt}\selectfont $1$} \end{array}\!\!\!\right)
            A^{m-1}B\bigg(\cos\alpha+\cos\Big(\frac{360^{\circ}}{n}-\alpha\Big)+\cdots+  \\[0.2cm]
    &\qquad\qquad\qquad\qquad +\cos\Big((n-1)\,\frac{360^{\circ}}{n}-\alpha\Big)\bigg)+ \\[0.2cm]
    &\qquad +\left(\!\!\!\begin{array}{c} \text{\fontsize{13}{14pt}\selectfont $m$} \\ \text{\fontsize{13}{14pt}\selectfont $2$} \end{array}\!\!\!\right)
            A^{m-2}B^2\bigg(\cos^2\alpha+\cos^2\Big(\frac{360^{\circ}}{n}-\alpha\Big)+\cdots+  \\[0.2cm]
    &\qquad\qquad\qquad\qquad +\cos^2\Big((n-1)\,\frac{360^{\circ}}{n}-\alpha\Big)\bigg)- \\[1cm]
    &\qquad -\left(\!\!\!\begin{array}{c} \text{\fontsize{13}{14pt}\selectfont $m$} \\ \text{\fontsize{13}{14pt}\selectfont $3$} \end{array}\!\!\!\right)
            A^{m-3}B^3\bigg(\cos^3\alpha+\cos^3\Big(\frac{360^{\circ}}{n}-\alpha\Big)+\cdots+  \\[0.2cm]
    &\qquad\qquad\qquad\qquad +\cos^3\Big((n-1)\,\frac{360^{\circ}}{n}-\alpha\Big)\bigg)+ \\[0.2cm]
    &\qquad\qquad\qquad \vdots \\[0.2cm]
    &\qquad \pm\left(\!\!\!\begin{array}{c} \text{\fontsize{13}{14pt}\selectfont $m$} \\ \text{\fontsize{13}{14pt}\selectfont $m$} \end{array}\!\!\!\right)
            B^m\bigg(\cos^m\alpha+\cos^m\Big(\frac{360^{\circ}}{n}-\alpha\Big)+\cdots+  \\[0.2cm]
    &\qquad\qquad\qquad\qquad +\cos^m\Big((n-1)\,\frac{360^{\circ}}{n}-\alpha\Big)\bigg).
\end{align*}
 }
\vskip+0.2cm
\noindent
By using the Lemma \ref{lem:2} each sum with negative sign ``{\fontsize{14}{16pt}\selectfont $-$}'' is zero, because they contain odd powers and only sum with even powers remains.

If {\fontsize{14}{16pt}\selectfont $m$} is even the sum is:
\vskip+0.02cm
{\fontsize{14}{16pt}\selectfont
\allowdisplaybreaks
\begin{align*}
    S_n^{(2m)} & =nA^m+  \\[0.2cm]
    &\qquad +\left(\!\!\!\begin{array}{c} \text{\fontsize{13}{14pt}\selectfont $m$} \\ \text{\fontsize{13}{14pt}\selectfont $2$} \end{array}\!\!\!\right)
            A^{m-2}B^2\bigg(\cos^2\alpha+\cos^2\Big(\frac{360^{\circ}}{n}-\alpha\Big)+\cdots+ \\[0.2cm]
    &\qquad\qquad\qquad\qquad +\cos^2\Big((n-1)\,\frac{360^{\circ}}{n}-\alpha\Big)\bigg)+ \\[0.2cm]
    &\qquad\qquad\qquad \vdots \\[0.2cm]
    &\qquad +\left(\!\!\!\begin{array}{c} \text{\fontsize{13}{14pt}\selectfont $m$} \\ \text{\fontsize{13}{14pt}\selectfont $m$} \end{array}\!\!\!\right)
            B^m\bigg(\cos^m\alpha+\cos^m\Big(\frac{360^{\circ}}{n}-\alpha\Big)+\cdots+ \\[0.2cm]
    &\qquad\qquad\qquad\qquad +\cos^m\Big((n-1)\,\frac{360^{\circ}}{n}-\alpha\Big)\bigg)= \\[0.2cm]
    & =n\Bigg(A^m+\sum_{k=1}^{\frac{m}{2}}
        \left(\!\!\!\begin{array}{c} \text{\fontsize{13}{14pt}\selectfont $m$} \\ \text{\fontsize{13}{14pt}\selectfont $2k$} \end{array}\!\!\!\right)
             A^{m-2k}B^{2k}\,\frac{1}{2^{2k}}
        \left(\!\!\!\begin{array}{c} \text{\fontsize{13}{14pt}\selectfont $2k$} \\ \text{\fontsize{13}{14pt}\selectfont $k$} \end{array}\!\!\!\right)\Bigg).
\end{align*}
 }
\vskip+0.2cm
\noindent
If {\fontsize{14}{16pt}\selectfont $m$} is odd the sum is
\vskip+0.02cm
{\fontsize{14}{16pt}\selectfont
\allowdisplaybreaks
\begin{align*}
    S_n^{(2m)} & =nA^m+ \\[0.2cm]
    &\quad +\left(\!\!\!\begin{array}{c} \text{\fontsize{13}{14pt}\selectfont $m$} \\ \text{\fontsize{13}{14pt}\selectfont $2$} \end{array}\!\!\!\right)
            A^{m-2}B^2\bigg(\cos^2\alpha+\cos^2\Big(\frac{360^{\circ}}{n}-\alpha\Big)+\cdots+  \\[0.4cm]
    &\qquad\qquad\qquad\qquad +\cos^2\Big((n-1)\,\frac{360^{\circ}}{n}-\alpha\Big)\bigg)+\cdots+ \\[0.4cm]
    &\quad +\left(\!\!\!\begin{array}{c} \text{\fontsize{13}{14pt}\selectfont $m$} \\ \text{\fontsize{13}{14pt}\selectfont $m-1$} \end{array}\!\!\!\right)
            AB^{m-1}\bigg(\cos^{m-1}\alpha+\cos^{m-1}\Big(\frac{360^{\circ}}{n}-\alpha\Big)+ \\[0.4cm]
    &\qquad\qquad\qquad +\cdots+\cos^{m-1}\Big((n-1)\,\frac{360^{\circ}}{n}-\alpha\Big)\bigg)= \\[0.4cm]
    & =n\Bigg(A^m+\sum_{k=1}^{\frac{m-1}{2}}
        \left(\!\!\!\begin{array}{c} \text{\fontsize{13}{14pt}\selectfont $m$} \\ \text{\fontsize{13}{14pt}\selectfont $2k$} \end{array}\!\!\!\right)
                     A^{m-2k}B^{2k}\,\frac{1}{2^{2k}}
        \left(\!\!\!\begin{array}{c} \text{\fontsize{13}{14pt}\selectfont $2k$} \\ \text{\fontsize{13}{14pt}\selectfont $k$} \end{array}\!\!\!\right)\Bigg).
\end{align*}
 }
\vskip+0.4cm

Using integer part, it is possible to write obtained results in one formula:
\vskip+0.04cm
{\fontsize{14}{16pt}\selectfont
$$  S_n^{(2m)}=n\Bigg(A^m+\sum_{k=1}^{[\frac{m}{2}]}
            \left(\!\!\!\begin{array}{c} \text{\fontsize{13}{14pt}\selectfont $m$} \\ \text{\fontsize{13}{14pt}\selectfont $2k$} \end{array}\!\!\!\right)
                        A^{m-2k}B^{2k}\,\frac{1}{2^{2k}}\,
            \left(\!\!\!\begin{array}{c} \text{\fontsize{13}{14pt}\selectfont $2k$} \\ \text{\fontsize{13}{14pt}\selectfont $k$} \end{array}\!\!\!\right)\Bigg),     $$
 }
\vskip+0.4cm
\noindent
which proves the Theorem \ref{th:8.1}.

\bigskip
\bigskip
\bigskip

Thus, for each even {\fontsize{14}{16pt}\selectfont $(2m)$}-th power {\fontsize{14}{16pt}\selectfont $m=1,2,\dots,n-1$}:
\vskip+0.04cm
{\fontsize{14}{16pt}\selectfont
\begin{multline*}
    \qquad \Big\{X:\;\;\sum_{i=1}^n d^{2m}(X,P_i)=S_n^{(2m)}=const\Big\}= \\
    =\begin{cases}
                    \text{\fontsize{12}{14pt}\selectfont circle}, & \text{\fontsize{12}{14pt}\selectfont if}\;\; S_n^{(2m)}>nr^{2m} \\[0.15cm]
                    {\fontsize{12}{14pt}\selectfont \text{point}}, & \text{\fontsize{12}{14pt}\selectfont if}\;\; S_n^{(2m)}=nr^{2m} \\[0.15cm]
                    \;\;\varnothing\;\;, & \text{\fontsize{12}{14pt}\selectfont if}\;\; S_n^{(2m)}<nr^{2m}
                \end{cases}\,.  \qquad
\end{multline*}
}
\vskip+0.4cm

In case of circle, the radius {\fontsize{14}{16pt}\selectfont $\ell$} satisfies:
\vskip+0.04cm
{\fontsize{14}{16pt}\selectfont
$$  S_n^{(2m)}=n\Bigg[(r^2+\ell^2)^m+\sum_{k=1}^{[\frac{m}{2}]}
            \left(\!\!\!\begin{array}{c} \text{\fontsize{13}{14pt}\selectfont $m$} \\ \text{\fontsize{13}{14pt}\selectfont $2k$} \end{array}\!\!\!\right)
                         (r^2+\ell^2)^{m-2k}(r\ell)^{2k}
            \left(\!\!\!\begin{array}{c} \text{\fontsize{13}{14pt}\selectfont $2k$} \\ \text{\fontsize{13}{14pt}\selectfont $k$} \end{array}\!\!\!\right)\Bigg].       $$
 }
\vskip+0.4cm

What happens if we consider the sum of even powers more than {\fontsize{14}{16pt}\selectfont $2(n-1)$}? In these cases, cosine functions appear in the sum expressions (see note of Lemma \ref{lem:1}), which proves that the locus is not a circle. But what kind of figures are there? The investigation of this subject is the aim of future studies.

\vskip+3cm
\section*{References}
\vskip+0.5cm

\begin{itemize}
\item[{[1]}] Prosolov Viktor,\textit{ Problems in Plane and Solid Geometry}. Vol. 1. \textit{Plane Geometry}. Translated and edit by Dimitry Leites, 2005.

\medskip

\item[{[2]}] Jones, Dustin. A Collection of Loci Using Two Fixed Points. \emph{Missouri Journal of Mathematical Sciences}, No. 19, 2007.

\end{itemize}

\vskip+3cm

\noindent \textbf{Author's address:}

\medskip

\noindent {Georgian-American High School, 18 Chkondideli Str., Tbilisi 0180, Georgia.}

\noindent {\small \textit{E-mail:} \texttt{mathmamuka@gmail.com} }

\end{document}